\documentclass[12pt,twoside]{amsart}
\usepackage{times,amssymb,hyperref,t1enc,a4wide}

\newtheorem{theorem}{Theorem}[section]
\newtheorem{definition}[theorem]{Definition}
\newtheorem{proposition}[theorem]{Proposition}

\newtheorem{lemma}[theorem]{Lemma}
\theoremstyle{remark}

\newcommand{\la}{\langle}
\newcommand{\ra}{\rangle}
\newcommand{\C}{\mathbb{C}}
\newcommand{\Z}{\mathcal {Z}}

\newcommand{\R}{\mathbb{R}}
\newcommand{\N}{\mathbb{N}}
\newcommand{\Sph}{\mathbb S}

\def \proof {\noindent {\bf Proof.}\ \ }
\def \qed {{\mbox{}\nolinebreak\hfill\rule{2mm}{2mm}\par\medbreak} }

\begin{document}
\title{Sufficient conditions for Sampling and interpolation on the sphere}
\author{Jordi Marzo}
\address{Dept.\ Matem\`atica Aplicada i An\`alisi,
 Universitat  de Barcelona,
Gran Via 585, 08007 Bar\-ce\-lo\-na, Spain}
\email{jmarzo@ub.edu}
\author{Bharti Pridhnani}
\address{Dept.\ Matem\`atica Aplicada i An\`alisi,
 Universitat  de Barcelona,
Gran Via 585, 08007 Bar\-ce\-lo\-na, Spain}
\email{bharti.pridhnani@ub.edu}


\thanks{}

\date{\today}
\begin{abstract}
     We obtain sufficient conditions for 
      arrays of points, $\mathcal{Z}=\{\mathcal{Z}(L) \}_{L\ge 1},$  on the unit sphere $\mathcal{Z}(L)\subset \mathbb{S}^d,$ to be
     Marcinkiewicz-Zygmund and interpolating arrays 
     for spaces of spherical harmonics. The conditions are in terms of the mesh norm and the separation radius of $\mathcal{Z}(L).$
\end{abstract}

\maketitle

\section{Introduction}

Let $L^{p}(\mathbb{S}^{d})$ the Banach space
of measurable functions defined in the unit sphere $\mathbb{S}^{d}$ with
\[
\| f \|_{p}^{p}=\int_{\mathbb{S}^{d}}|f(z)|^{p}d\sigma(z)<\infty,
\]
if $1\le p< \infty,$ and
\[
\| f \|_{\infty}=\sup_{z\in \mathbb{S}^{d}}|f(z)|<\infty,
\]
when $p=\infty$. Here $\sigma$ stands for the Lebesgue surface measure in $\mathbb{S}^d.$\\

For $l\geq 0$ an integer, let $\mathcal H_l$ be the space of spherical harmonics of degree $l$ in 
$\Sph^d$ i.e. the space of eigenfunctions of the 
Laplace-Beltrami operator on $\Sph^d$
\[
\Delta_{\Sph^d} Y+l(l+d-1)Y=0,\;\; Y\in \mathcal H_l.
\] 
We denote by $\Pi_L,$ for $L$ a nonnegative integer, the space of spherical harmonics of 
degree not exceeding $L$
\[
\Pi_L=\text{span}\bigcup^{L}_{l=0}\mathcal H_l.
\]
With respect to the inner product in $L^2(\Sph^d)$ the spaces $\mathcal H_l$ are orthogonal. We 
denote by $h_l$ and $\pi_L$ the dimension of $\mathcal H_l$ and $\Pi_L$, respectively. By 
Stirling's formula, $\pi_L\simeq L^d$, when $L\to \infty$. Let $Y^{1}_{l},\ldots,Y^{h_l}_{l}$ be 
an orthonormal basis of $\mathcal H_l.$ The reproducing kernel in $\Pi_L$ is given by
\[
K_L(u,v)=\sum^{L}_{l=0}\sum^{h_l}_{j=1}Y^{j}_{l}(u)\overline{Y^j_{l}(v)}=
C_{d,L} P^{(1+\lambda,\lambda)}_L(\langle u,v\rangle),\;\; u,v\in \Sph^d,
\]
where $\la u, v\ra$ stands for the scalar product in $\R^{d+1},$ $d=2\lambda+2,$ $C_{d,L}L^{-d/2}$ goes to a positive constant when $L\to+\infty$ and $P^{(\alpha,\beta)}_L$ 
are the Jacobi polynomials of degree $L$ and index $(\alpha,\beta),$ normalized so that 
\[
P^{(\alpha,\beta)}_L(1)=\left(\begin{array}{c} 
L+\alpha\\
L
\end{array}\right).
\]
  
To discretize the $L^p$-norms in the space of spherical harmonics, we consider arrays of points on the sphere. More precisely, for any degree $L$ we take 
$m_L$ points in $\Sph^d$
\[
\mathcal{Z}(L)=\{ z_{L j}\in \mathbb{S}^{d}: 1\le j\le m_{L}\}, \;\; L\ge 0.
\]
This yields an array of points $\Z=\left\{\Z(L)\right\}_{L\geq 0}$ in $\Sph^d$. 
Denote $d(u,v)=\arccos \la u, v\ra$ the geodesic distance between $u,v\in \mathbb{S}^{d}.$

\begin{definition}
    Let $\mathcal{Z}=\{ \mathcal{Z}(L) \}_{L\ge 0}$ be an array with $m_{L}\ge \pi_{L}$ for all $L.$
    We call $\mathcal{Z}$ an $L^{p}$-Marcinkiewicz-Zygmund array, denoted by $L^{p}$-MZ, if 
    there exists a constant $C_{p}>0$ such that for all $L\ge 0$ and $Q\in \Pi_{L},$
\begin{equation}\label{MZ-ineq}
    \frac{C_{p}^{-1}}{\dim \Pi_{L}}\sum_{j=1}^{m_{L}}|Q(z_{L j})|^{p}
    \le \int_{\mathbb{S}^{d}}|Q(\omega)|^{p}d\sigma(\omega)
    \le \frac{C_{p}}{\dim \Pi_{L}}\sum_{j=1}^{m_{L}}|Q(z_{Lj})|^{p},
\end{equation}
if $1\le p<\infty,$ and
$$\sup_{\omega \in \mathbb{S}^{d}}|Q(\omega)|\le C\sup_{j=1,\dots , m_{L}}|Q(z_{Lj})|,$$
when $p=\infty.$
\end{definition}

In other words, the $L^p$-norm in $\Sph^d$ of a polynomial of degree $L$ is comparable to the discrete version
given by the weighted $\ell^p$-norm of its restriction to $ \mathcal{Z}(L).$ 
For the unit circle, $d=1,$ the spherical harmonics are trigonometric polynomials. In this case, for $m_L=\pi_L=2L+1,$ 
J. Marcinkiewicz and A. Zygmund proved that the array of roots of unity form an $L^p$-MZ array, \cite{MZ37},
observe that $C_2=1.$ In higher dimensions the situation is more delicate. For $m_L=\pi_L$ there are no $L^p$-MZ arrays when $p\neq 2$
and the case $p=2$ is open, see \cite{Mar07}. For $m_L$ big enough there are always $L^p$-MZ arrays, see for example \cite{MNW00, FM11}.

\begin{definition}
Let $\mathcal{Z}=\{ \mathcal{Z}(L) \}_{L\ge 0}$ be a triangular array with $m_L\le \pi_L$ for all $L.$ 
We say that $\mathcal{Z}$ is $L^p$- interpolating if, for arrays $\{ c_{Lj} \}_{L\ge 0, 1\le j\le m_j}$ of 
complex values such that
$$\sup_{L\ge 0}\frac{1}{\pi_L}\sum_{j=1}^{m_L}|c_{Lj}|^p<\infty,$$
there exists a sequence of polynomials $Q_L \in \Pi_L$ uniformly bounded in $L^p$ such that 
$$Q(z_{Lj})=c_{Lj},\;\; 1\le j\le m_j.$$
\end{definition}

The concept of MZ and interpolating families is important in signal processing. Somehow these concepts are opposite in nature. 
MZ families are dense enough so that the $L^p$-norm in $\Sph^d$ is comparable to a discrete version. 
On the other hand, interpolating families 
are sparse enough so that one can interpolate some given data. 

Observe that $\Z$ is $L^2$-MZ if and only if 
the normalized reproducing kernels of $\Pi_L$ at the points $\Z(L)$ form a frame with frame bounds 
independent of $L.$ Therefore, $\Z$ is $L^2$-MZ when $\Z(L)$ is a set of sampling for $\Pi_L$ with constants independent of $L.$  
Similarly, $\Z$ is $L^2$-interpolating if and only if the normalized reproducing kernel of $\Pi_L$, $k_L,$
at the points $\Z(L)$
form a Riesz sequence i.e.
\begin{equation}							\label{rieszseq}
C^{-1} \sum_{j=1}^{m_L} |a_{Lj}|^2 \le \int_{\mathbb{S}^d} \left| \sum_{j=1}^{m_L} a_{Lj} k_L(z, z_{L,j})\right|^2 d\sigma(z)\le C \sum_{j=1}^{m_L} |a_{Lj}|^2,
\end{equation}
for any $\{ a_{Lj} \}_{L,j}$
with $C>0$ independent of $L.$ 
When $\Z$ is both $L^2$-interpolating and $L^2$-MZ 
the normalized reproducing kernels at the points $\Z(L)$ form a Riesz basis, for more about these concepts see \cite{Sei95}.

We denote by $d(u,v)=\arccos \la u, v\ra$ the geodesic distance between $u,v\in \mathbb{S}^{d}.$

\begin{definition}
    An array $\mathcal{Z}=\{ \mathcal{Z}(L) \}_{L\ge 0}$ is uniformly separated if there is a positive number $\epsilon>0$
    such that
    $$d(z_{Lj},z_{Lk})\ge
    \frac{\epsilon}{L+1},\;\; \mbox{if}\;\; j\neq k,$$
    for all $L\ge 0.$
\end{definition}

The left hand side inequality in (\ref{MZ-ineq}) holds if and only if $\mathcal{Z}$ is a finite 
union of uniformly separated arrays, 
also the $L^p$ version of the right hand side inequality in (\ref{rieszseq}) holds if and only if $\mathcal{Z}$ is uniformly separated,
see \cite{Mar07} or \cite{OCP11} 
for the general case of a compact Riemannian manifold.

\begin{definition}
    Let $X$ be a subset of $\mathbb{S}^d.$ The mesh norm of $X$ is
$$\rho (X)=\sup_{u\in \mathbb{S}^{d}} d(u,X)=\sup_{u\in \mathbb{S}^{d}}  \inf_{z\in X} d(u,z).$$
    The separation radius of $X$ is
$$\delta(X)=\inf_{u\in X}  \inf_{v\in X\setminus\{ u \}} d(u,v).$$
\end{definition}

The mesh norm of $X\subset \mathbb{S}^d$ is therefore the maximal radius of a spherical cap which does not contain points from $X,$ and the 
separation radius is the minimal distance between points in $X.$ 

\subsection{Main results}

  Our results in this paper are the following sufficient conditions for an array $\Z$ to be $L^p$-MZ or $L^p$-interpolating.

\begin{theorem}\label{main_theo2}
Let $1\le p\le \infty$ and $\mathcal{Z}=\{ \mathcal{Z}(L) \}_{L\ge 0}$ be an array in $\mathbb{S}^d$ such that for all $L\ge 0$
\[
\delta(\Z(L))>\frac{\theta}{L},
\]
where $\theta>2j_\lambda,$ $j_\lambda$ is the first zero of the Bessel function $J_\lambda(t),$ and $d=2\lambda+2.$ Then
$\mathcal{Z}$ is $L^p$-interpolating family.
\end{theorem}

\begin{theorem}\label{main_theo}
Let $1\le p\le \infty$ and $\mathcal{Z}=\{ \mathcal{Z}(L) \}_{L\ge 0}$ be a uniformly separated array in $\mathbb{S}^d$ such that for all $L\ge 0$
\begin{equation}
    \rho(\mathcal{Z}(L))<\frac{\pi}{2L},
\end{equation}
    then $\mathcal{Z}$ is an $L^{p}$-Marcinkiewicz-Zygmund array.
\end{theorem}

We want to point out that there are other results about sufficient conditions for $L^p$-MZ and interpolating families of points 
on compact manifolds, but such conditions do not provide precise constants, see \cite{FM11,FM10,MNSW02,OCP12}.
We observe that due to the result mentioned above about minimal $L^p$-MZ (or maximal $L^p$-interpolating) arrays, an array
with $m_L=\pi_L$ cannot satisfy the conditions of Theorems \ref{main_theo},\ref{main_theo2}. 

When $\Sph^2$ and for some particular arrays of points, there are some results
about separation radius and mesh norm \cite{Rei90,SW04,DM05}. The results we know are not very precise and 
we would just mention one to illustrate the use of our results.
 
  The set $X=\{x_1,\dots , x_N \}\subset \Sph^2$
is said to be in $s$-extremal configuration if $X$ maximizes the Riesz $s$-energy
$$E(X)=\sum_{1\le i<j\le N}\frac{1}{|x_i-x_j|^s}$$
for subsets of $N$ points on the sphere. In \cite{KS98} it is assumed in order to get an estimate for the separation radius that te Voronoi cells around points 
in $s$-extremal configuration are all hexagons. Then it is obtained that
$$\delta(X_N)\sim \left( \frac{8\pi}{\sqrt{3}}\right)^{1/2}N^{-1/2}.$$ The same way it seems reasonable to estimate the mesh norm by the value
of the maximal radius of the hexagon getting
$$\rho(X_N)\sim \frac{1}{\sqrt{3}}\left( \frac{8\pi}{\sqrt{3}}\right)^{1/2} N^{-1/2}$$
  Therefore, if we take $(k L)^2$ points for degree $L$ in order to assure we get an $L^p$-interpolating array we need
  $k<0.792,$ and to get an $L^p$-MZ array we need $k>1.4.$


\subsection{Outline of the paper}

In Section 2 we prove Theorem \ref{main_theo2}. We use the classical approach by Ingham to obtain sufficient conditions for interpolation, \cite{Ing36}.
Ingham idea has been used in different context, \cite{OU10,KL05}.
The main problem is the construction of appropiate pick functions, Lemma \ref{pickfunc}, which depend on 
an upper bound for the first eigenvalue of a spherical cap, \cite{Pin81,BCG83}.
Estimates for the first eigenvalue are known also in general Riemannian manifols. 

In Section 3 we prove Theorem \ref{main_theo}. Our approach
follow the classical ideas of Beurling to study sampling sequences in Bernstein space, \cite{Beu89}. 
We define weak limits of an array, and relate uniqueness sets
with the $L^\infty$-MZ property. Our result is consequence of a uniqueness result due also to Beurling. 

In what follows, when we write $A\lesssim B$, $A\gtrsim B$ or $A\simeq B$ , we mean that there are constants independent 
of $L$ such that $A\leq CB$, $A\geq CB$ or $C_1B\leq A\leq C_2B$, respectively. Also, the value of the constants 
appearing during a proof may change but they will be still denoted with the same letter.


\section{Sufficient condition for interpolation}

In this section we prove Theorem \ref{main_theo2}. We adapt a nice idea of Ingham for Dirichlet series, \cite{Ing36}.
The idea rely on the construction of some pick functions with appropiate spectral properties. The existence of such functions can be established 
studying the first eigenvalue/eigenfunction of the Laplace-Beltrami operator on the sphere 
from which a spherical cap has been removed, see \cite{KL05,OU10}.
Estimates on this first eigenvalue, \cite{FH76, Pin81, BCG83}, together with a result about perturbation of interpolating 
arrays, \cite[Lemma 4.11]{Mar07}, provide the result. For the unit circle, $d=1,$ the condition in Theorem \ref{main_theo2} is $\theta=\pi,$
and the proof is technically simpler.

To prove Theorem \ref{main_theo2} we use the functions given by the following lemma. 

\begin{lemma}											\label{pickfunc}
Given $L$ and $\theta$ as in Theorem \ref{main_theo2}, there exists functions $F_L$ such that
\begin{enumerate}
	\item $F_{L}\in L^{2}([-1,1]).$
	\item $ \mbox{\rm supp}\, F_L\subset [\cos(\theta/L),1]$.
	\item $[F_L (\langle u,\cdot\rangle )](\ell,j)\le 0$ for all $\ell>L$ and 
	$[F_L (\langle u,\cdot\rangle )](\ell,j)\lesssim \theta^{d/2}$ for all $\ell\leq L$.
	\item $F_L(1)\simeq \pi_L$.
\end{enumerate}
  where $[F_L (\langle u,\cdot\rangle )](\ell,j)$ stands for the Fourier coefficient $\int_{\Sph^d} F_L(\langle z_{Lj},u\rangle) \overline{Y_{\ell}^j(u)}d\sigma(u).$
\end{lemma}

Before establishing the existence of such functions we prove our main result.
\vskip 0.2cm


\proof [Theorem \ref{main_theo2}]
Recall that the normalized reproducing kernel can be written as
\[
k_L(z, z_{L,j})=\frac{1}{\sqrt{\pi_L}}\sum_{\ell=0}^L \sum_{k=1}^{h_\ell}Y_{\ell}^k(z)\overline{Y_{\ell}^k(z_{Lj})}.
\]
Let $F_L(x)$ be given by Lemma \ref{pickfunc}, i.e. $F_L$ is a continuous function defined in $-1\le x \le 1$ such that
\[
\mbox{supp}\; F_L\subset [\cos (\theta/L),1],
\]
and $F_L(1)\ge \pi_L \sim L^d$. Moreover, the Fourier coefficients of $F_L$ are negative for $\ell>L$ and uniformly bounded 
by $C\theta^{d/2}$ for $\ell\leq L$. Thus by using these estimates and Funk-Hecke formula we get 
\[
\int_{\mathbb{S}^d} \left| \sum_{j=1}^{m_L} c_{L,j} k_L(z, z_{L,j})\right|^2 d\sigma(z)=
\frac{1}{\pi_L}\int_{\mathbb{S}^d}\left| \sum_{\ell=0}^L \sum_{k=1}^{h_\ell} \left( \sum_{j=1}^{m_L} c_{Lj} \overline{Y_{\ell}^k(z_{Lj})}\right)
Y_{\ell}^k(z) \right|^2 d\sigma(z)
\]
\[
=\frac{1}{\pi_L} \sum_{\ell=0}^L \sum_{k=1}^{h_\ell} \left| \sum_{j=1}^{m_L} c_{Lj} \overline{Y_{\ell}^k(z_{Lj})} \right|^2
\]
\[
\gtrsim \frac{1}{\theta^{d/2}}\frac{1}{\pi_L} \sum_{\ell=0}^{+\infty}  \sum_{k=1}^{h_\ell} \sum_{i,j=1}^{m_L} [F_L(\langle z_{Lj},\cdot\rangle)](\ell,k)
 c_{Lj}\overline{c_{Li}} \overline{Y_{\ell}^k(z_{Lj})}Y_{\ell}^k(z_{Li})
 \]
 \[
=\frac{1}{\theta^{d/2}}\frac{1}{\pi_L}\sum_{i,j=1}^{m_L}  c_{Lj}\overline{c_{Li}} \sum_{\ell=0}^{+\infty} \sum_{k=1}^{h_\ell} 
[F_L(\langle z_{Lj},\cdot\rangle)](\ell,k)  \overline{Y_{\ell}^k(z_{Lj})}Y_{\ell}^k(z_{Li})
\]
\[
=\frac{1}{\theta^{d/2}}\frac{1}{\pi_L}\sum_{i,j=1}^{m_L}  c_{Lj}\overline{c_{Li}} F_L(\cos d(z_{Li},z_{Lj}))=
\frac{1}{\theta^{d/2}}\frac{F_L(1)}{\pi_L}\sum_{j=1}^{m_L}  |c_{Lj}|^2\gtrsim  \sum_{j=1}^{m_L}  |c_{Lj}|^2,
\]
and $\Z$ is therefore $L^2$-interpolating, because the other inequality in (\ref{rieszseq}) follows directly from the separation.

In order to prove the result for other $p\neq 2$ we define, 
for $\delta>0$, the per	turbed array $\Z_{\delta}=\left\{\Z(L_{1+\delta})\right\}_{L}$, where $L_{1+\delta}=[L(1+\delta)]$. It was 
proved in \cite[Lemma 4.11]{Mar07} that if $\Z$ is $L^2$-interpolating then $\Z_{-\delta}$ is $L^p$-interpolating for all $p\in[1,\infty]$. Therefore, 
assume that $\Z$ satisfies the geometric separation condition
\[
\eta:=L\min_{i\neq j} d(z_{Li},z_{Lj})>\theta.
\]
Let $\delta>0$ be small enough so that $\eta>\theta+\delta \theta$. We assume that $L>>1$ so that $L\delta >1$. Then
\[
Ld(z_{L_{1+\delta} i}, z_{L_{1+\delta} j})>\frac{L}{L_{1+\delta}}\eta>\theta \frac{L}{L_{1+\delta}}(1+\delta)\geq \theta.
\]
Thus, the perturbed array $\Z_{\delta}$ satisfies the same separation condition as $\Z,$ then $\Z_{\delta}$ is $L^2$-interpolating and $\Z=(\Z_\delta)_{-\delta}$ is 
$L^p$-interpolating for all $1\leq p\leq \infty$. 
\qed


\proof [Lemma \ref{pickfunc}]
Let $C_{\theta/2L}$ be the spherical cap of those $x=(x_1,\dots , x_{d+1})\in \Sph^d$ such that
$\cos\frac{\theta}{2L}<x_{d+1}\leq 1$. We denote by $f_0$ the eigenfunction of the problem
\[
\Delta_{\Sph^d} f_0+\lambda_{0,L}f_0=0
\]
in $C_{\theta/2L}$ corresponding to the first eigenvalue $\lambda_{0,L}$ of the Laplace-Beltrami operator $\Delta_{\Sph^d}$. It is known, see \cite{FH76}, 
that $f_0$ belongs to the class of zonal functions, Lipschitzian, nonnegative, non identically zero and with support in $[0,\frac{\theta}{2L}].$	 

We normalize $\|f_0\|^2_2\simeq \pi_L$ and define the zonal function
\[
F_L=\left(1+\frac{\Delta_{\Sph^d}}{L(L+d-1)}\right)(f_0\ast f_0).	
\]
with support in $[\cos \frac{\theta}{L},1].$

Let $\{ Y_\ell^j \}_{j,\ell}$ be the orthonormal basis in $L^2(\Sph^d)$ given by the spherical harmonics, then
$$F_L(\langle u,v\rangle)=\sum_{\ell\geq 0}\sum^{h_\ell}_{j=1}  [F_L(\langle u, \cdot\rangle)] (\ell,j)Y^{j}_{\ell}(v),$$
where by Funk-Hecke
\[
[F_L(\langle u, \cdot\rangle)] (\ell,j)=\int_{\Sph^d}F_L(\langle u,v \rangle)\overline{Y^{j}_{l}(v)}d\sigma(v)=\widehat{F_L}(\ell)\overline{Y^{j}_{l}(u)},
\]
and
\[
\widehat{F_L}(\ell)=\frac{\sigma(\Sph^d)}{C^{(d-1)/2}_{\ell}(1)}\frac{\int^{1}_{-1}F_L (t)C^{(d-1)/2}_{\ell}(t)(1-t^2)^{(d-2)/2}dt}
{\int^{1}_{-1}(1-t^2)^{(d-2)/2}dt},
\]
here $C^{\alpha}_\ell$ is the Gegenbauer polynomial of order $\alpha$ and degree $\ell$.

On the other hand
\begin{align*}
[F_L(\langle u, \cdot\rangle)] (\ell,j)& =\int_{\Sph^d}\left(1+\frac{\Delta_{\Sph^d}}{L(L+d-1)}\right)(f_0\ast f_0)(\langle u,v\rangle)\overline{Y^{j}_{\ell}(v)}d\sigma(v)\\
&=\int_{\Sph^d}(f_0\ast f_0)(\langle u,v\rangle)\overline{Y^{j}_{\ell}(v)}d\sigma(v)+
\int_{\Sph^d}\frac{\Delta_{\Sph^d}(f_0\ast f_0)(\langle u,v\rangle)}{L(L+d-1)}\overline{Y^{j}_{\ell}(v)}d\sigma(v)\\
&=\widehat{f_0}(\ell)^2\overline{Y^{j}_{\ell}(u)}
+\frac{1}{L(L+d-1)}\overline{\int_{\Sph^d}(f_0\ast f_0)(\langle u,v\rangle)\Delta_{\Sph^d} Y^{j}_{\ell}(v) d\sigma(v)}\\
&=\widehat{f_0}(\ell)^2\overline{Y^{j}_{\ell}(u)}
-\frac{\ell(\ell+d-1)}{L(L+d-1)}\int_{\Sph^d}(f_0\ast f_0)(\langle u,v\rangle)\overline{Y^{j}_{\ell}(v)}d\sigma(v)\\
&=\left(1-\frac{\ell(\ell+d-1)}{L(L+d-1)}\right)\widehat{f_0}(\ell)^2\overline{Y^{j}_{\ell}(u)}.
\end{align*}

So we have proved that
$$\widehat{F_L}(\ell)=\left(1-\frac{\ell(\ell+d-1)}{L(L+d-1)}\right)\widehat{f_0}(\ell)^2$$

Note that the coefficients $\widehat{F_L}(\ell) \leq 0$ for all $\ell>L$. Now we are going to prove that the coefficients 
$\widehat{f_0}(\ell)$ are bounded for $\ell \leq L$. 
\begin{align*}
|\widehat{f_0}(\ell)|& \simeq \frac{1}{C^{(d-1)/2}_{\ell}(1)}\left|\int^{1}_{-1}f_0(t)C^{(d-1)/2}_{\ell}(t)(1-t^2)^{(d-2)/2}dt\right|\\
& \lesssim \frac{1}{C^{(d-1)/2}_{\ell}(1)}\|f_0\|_2\left(\int^{\theta/2L}_{0}C^{(d-1)/2}_{\ell}(\cos\theta)^2\sin^{d-1}\theta d\theta\right)^{1/2}\\
& \simeq \frac{\sqrt{\pi_L}}{C^{(d-1)/2}_{\ell}(1)}\left(\int^{\theta/2L}_{0}C^{(d-1)/2}_{\ell}(\cos\theta)^2\sin^{d-1}\theta d\theta\right)^{1/2}\lesssim \theta^{d/2},
\end{align*}
where we have used that $C^{(d-1)/2}_{\ell}(1)\simeq \ell^{(d-2)/2}$ and $C^{(d-1)/2}_{\ell}(t)\leq C\ell^{\max(\frac{d-2}{2},-\frac{1}{2})}$ 
(see \cite[Section 7.32]{Sze39}).

On the other hand, note that
$$(f_0 \ast f_0)(1)=(f_0\ast f_0)(\langle N,N\rangle)=\int^{\pi}_{0}f^2_0(\cos\theta)\sin^{d-1}\theta d\theta\simeq \pi_L.$$
Thus,
\begin{align*}
F_L(1)& =(f_0\ast f_0)(1)+\frac{1}{L(L+d-1)}\Delta_{\Sph^d}(f_0\ast f_0)(1)\\
&=(f_0\ast f_0)(1)+\frac{1}{L(L+d-1)}(f_0\ast \Delta_{\Sph^d} f_0)(1)\\
&=\left(1-\frac{\lambda_{0,L}}{L(L+d-1)}\right)(f_0\ast f_0)(1)\simeq \pi_L\left(1-\frac{\lambda_{0,L}}{L(L+d-1)}\right).
\end{align*}
So we need to find the smallest $\theta$ so that the quantity
\[
1-\frac{\lambda_{0,L}}{L(L+d-1)}>0.
\]
Equivalently, we need the smallest $\theta$ so that $\lambda_{0,L}<L(L+d-1)$. Using the upper bound from \cite{BCG83}, we get that
\[
\lambda_{0,L} \frac{\theta^2}{4L^2} <j^2_{\frac{d-2}{2}},
\]
where $j_{(d-2)/2}$ is the first zero of the Bessel function $J_{(d-2)/2}$. So taking $\theta$ as in the hypothesis we have the result. 
\qed


\section{Sufficient condition for Sampling}

In this section we follow the classical approach used by Beurling to study sampling sequences in the space of bounded 
bandlimited functions, \cite{Beu89}. First we identify the space of spherical harmonics composed with the exponential map as a subspace of 
the space of bounded bandlimited functions. Then we define the concept of weak limit of an array, and relate uniqueness sets
with $L^\infty$-MZ arrays. Finally, we get
Theorem \ref{main_theo} by using a result of Beurling about uniqueness sets and a result about perturbation of MZ arrays, \cite[Lemma 4.9.]{Mar07}.

Let $x=(x_{1},\dots ,x_d)\in \R^d$ and denote $\phi(x)=(x_{1}^{2}+\dots +x_{d}^{2})^{1/2}.$
The exponential map in $\Sph^d$ is defined by
$$\exp(x)=\left( x_1 \frac{\sin \phi(x)}{\phi(x)},\dots , x_d \frac{\sin \phi(x)}{\phi(x)}, \cos \phi(x) \right)\in \Sph^d.$$
Observe that $\exp(z)$ is defined also for $z\in \C^d$ and is an entire function.

Given $Q_{L}\in\Pi_{L}$ we define the function
\begin{equation}\label{funcions}
    \widetilde{Q}_{L}(z)=Q_{L}(\exp (z/L)),\;\; z\in \C^{d},
\end{equation}
and the corresponding space $\widetilde{\Pi}_{L}.$ Observe that
$$\sup_{u\in \mathbb{S}^{d}}|Q_{L}(u)|=\sup_{x\in \R^{d}}|\widetilde{Q}_{L}(x)|.$$

    The following result shows that
    functions in $\widetilde{\Pi}_{L}$  are entire in $\C^{d}$
    with Fourier-Laplace transform supported in the unit ball.

\begin{proposition}
    If $Q\in \Pi_{L}$ then the Fourier transform of $\widetilde{Q}$ has support in the unit ball of $\R^d.$
\end{proposition}

\proof
    The reproducing kernel of $\Pi_{L}$ centered at $v\in \mathbb{S}^d$ is, up to constants,
    the Jacobi polynomial $u\mapsto P^{(1+\lambda,\lambda)}_{L}(\la u, v \ra).$
    Let $y\in \R^{d}$ be such that $\exp(y/L)=v.$
    Consider the entire function
    $$\C^{d} \ni z\mapsto P^{(1+\lambda,\lambda)}_{L}\left( \la \exp(z/L),\exp(y/L) \ra \right).$$
    It is enough to see that for
    $\ell \le L$ and for some constants $C,N\ge 0$ (that may depend on $L$)
    $$\left| \la \exp(z/L),\exp(y/L) \ra \right|^{\ell}\le C (1+|z|)^{N} e^{|\Im z|}, \;\;z\in \C^{d}.$$
    For any $\zeta \in \C$ one has $2 ( \Im \zeta )^2=|\zeta|^2-\Re \zeta^2 .$ Therefore for any $z=(z_{1},\dots , z_d)\in \C^{d}$
    $$\sum_{i=1}^d \Re z_{i}^2 + 2 \sum_{i=1}^d (\Im z_{i})^2=\sum_{i=1}^d |z_{i}|^2,$$
    and we get by the triangle inequality that
    $$|\Im (z_{1}^{2}+\dots +z_{d}^{2})^{1/2}|\le \left( \sum_{i=1}^{d}(\Im z_{i})^{2} \right)^{1/2}.$$
    Finally
\begin{align*}
    |\la \exp(z/L),\exp(y/L) \ra| & =\left| \frac{\la z, y\ra}{L|y|}
    \frac{\sin \phi(z)/L}{\phi(z)/L}\sin \frac{|y|}{L}+\cos \frac{\phi(z)}{L}\cos \frac{|y|}{L}\right|
    \\
    &
    \le
    \left(1+\left|
    \frac{\la z, y\ra}{L|y|}
    \sin \frac{|y|}{L}\right|\right)
    e^{|\Im \phi(z)/L|}\le
    C_{y,L}(1+|z|)e^{\phi(\Im z)/L},
\end{align*}
    and $\widetilde{Q}$ is the Fourier-Laplace transform of a distribution supported in the unit ball of $\R^d.$
\qed

    Let $\mathcal{B}$ be the Bernstein space of entire functions in $\C^{d},$ bounded in $\R^{d}$ with
    Fourier transform supported in the unit ball of $\R^d,$ endowed with the uniform norm.

    Given an array $\mathcal{Z}$ we send the points in $\mathcal{Z}(L)$ to $\R^{d}$ via the exponential map and define the
    corresponding family
    of weak limits.
    The
    Fr\'echet distance between the closed sets $A,B\subset \R^{d}$ is given by
    $$[A,B]=\inf_{t>0} \{A\subset B+B(0,t),B \subset A+B(0,t)\}.$$

\begin{definition}
    Let $\mathcal{Z}=\{ \mathcal{Z}(L) \}_{L\ge 0}$ be an
    array in $\mathbb{S}^d.$
    We say that $\Lambda\subset \R^{d}$ is a weak limit of
    $\mathcal{Z},$ denoted as $\Lambda\in W(\mathcal{Z}),$
    if there exist rotations $\rho_{L}\in SO(d+1),$ such that
    $$L\exp^{-1}(\rho_{L}\mathcal{Z}(L))\rightharpoonup \Lambda,$$
    where the above expression means that for any $K\subset \R^{d}$ compact
    $$[( L \exp^{-1} ( \rho_{L}\mathcal{Z}(L))\cap K)
    \cup \partial K,(\Lambda\cap K)\cup \partial K]\to 0,\;\;\; L\to \infty.$$
\end{definition}

\begin{proposition}                                                             \label{lim_feb}
    If any $\Lambda\in W(\mathcal{Z})$ is a uniqueness set for $\mathcal{B}$ then
    $\mathcal{Z}$ is a $L^{\infty}$-MZ array.
\end{proposition}

\proof
    We argue by contradiction. Suppose that
    $\Lambda\in W(\mathcal{Z})$ is a uniqueness set for
    $\mathcal{B}$ but $\mathcal{Z}$ is not $L^{\infty}$-MZ.
    For any $n\in \N$ there exists $Q_{n}\in \Pi_{L_{n}}$ such that $Q_{n}(N)=\|  Q_{n} \|_{\infty}=1$ and
\begin{equation}                                                                                    \label{contra1}
    \frac{1}{n}>\sup_{j=1,\dots ,m_{L_{n}}}|Q_{n}(z_{L_{n}j})|.
\end{equation}
    From  the sequence $( \widetilde{Q}_{n} )_{n}$ defined as in (\ref{funcions}) it is possible
    to select a subsequence (see \cite[3.3.6.]{Nik75}) converging uniformly on compact sets of $\C^{d}$ to some function
    $f\in \mathcal{B}$ with $f(N)=1.$ We denote this subsequence as before.
    For any $\lambda\in \Lambda$ there exists a sequence
    $z_{L_{n_{k}}j_{n_{k}}}\in \mathcal{Z}(L_{n_{k}})$ such that
    $$\R^{d}\ni w_{L_{n_{k}}}=L_{n_{k}}\exp^{-1}(z_{L_{n_{k}}j_{n_{k}}})\to \lambda,\;\;\; k\to \infty.$$
    We denote this subsequence as before and we get
    $$|f(\lambda)|\le |f(\lambda)-\widetilde{Q}_{n}(\lambda)|+|\widetilde{Q}_{n}(\lambda)-\widetilde{Q}_{n}(w_{L_{n}j_{n}})|+
    |\widetilde{Q}_{n}(w_{L_{n}j_{n}})|.$$
    The first term on the right side clearly goes to zero. Also the last term goes to zero because of (\ref{contra1}).
    For the second term, using Bernstein's inequality we get
\begin{align*}
    |\widetilde{Q}_{n}(\lambda)-\widetilde{Q}_{n}(\widetilde{z}_{L_{n}j_{n}})| & =
    |Q_{n}(\exp(\lambda/L_{n}))-Q_{n}(z_{L_{n}j_{n}})|
    \\
    &
    \le
    L_{n}d(\exp(\lambda/L_{n}),z_{L_{n}j_{n}})\| Q_{n} \|_{\infty}\to 0,\;\;\; n\to \infty,
\end{align*}
    because, see \cite[p. 229]{BC73},
    $$d(\exp \frac{x}{L},\exp \frac{y}{L})=\frac{|x-y|}{L}+o(L^{-1}),\;\;\;\; x,y\in \R^{d}.$$
    We get that $f$ vanish in $\Lambda,$ but $f\in \mathcal{B}$ and therefore $f= 0.$
\qed

    To prove our main result we use the following
    result about uniqueness due to Beurling, \cite[p. 310]{Beu89}.

\begin{theorem}[Beurling]                                  
    Let $f$ be an entire function in $\C^{d}.$ Assume that
    $$\limsup_{|\xi|\to \infty}\frac{\log |f(\xi)|}{|\xi|}=r<\infty,\;\;
    \int_{1}^{\infty}\max_{|\xi|\le t,\xi\in \R^{d}}\log |f(\xi)|\frac{dt}{t^{2}}<\infty.$$
    Assume $f=0$ on a discrete set $\Lambda\subset \R^{d}$ such that
    $$r\limsup_{\R^{d}\ni x \to \infty}\inf_{\lambda\in \Lambda}|x-\lambda|<\frac{\pi}{2}.$$
    Then $f= 0.$
\end{theorem}

\proof[Theorem \ref{main_theo}]
    Let $\mathcal{Z}$ be such that for all $L$ big enough
$$\rho(\mathcal{Z}(L))<\frac{\pi}{2L}.$$
    Let $\Lambda \in W(\mathcal{Z})$ and
    $X_{L}=L\exp^{-1}(\rho_{L}\mathcal{Z}(L))\subset \R^{d}$ with
    $X_L\rightharpoonup \Lambda.$
    We want to see that $\Lambda$ is a uniqueness set for $\mathcal{B}.$

    For $f\in \mathcal{B}$ and $\epsilon>0$ there exists $A_{\epsilon}>0$ such that
    $$|f(\xi)|\le A_{\epsilon}e^{(1+\epsilon)|\xi|},\;\; \xi \in \C^d.$$
    Also $|f(x)|\le M<\infty$ for $x\in \R^{d},$ so we can apply
    Beurling's result (with $r=1$) getting that any $\Lambda\subset \R^d$ such that
\begin{equation}                        \label{metric}
    \limsup_{\R^{d}\ni x \to \infty}\inf_{\lambda\in \Lambda}|x-\lambda|<\frac{\pi}{2},
\end{equation}
     is a uniqueness set for $\mathcal{B}.$

    If
    $$\sup_{|x|<\pi L}     \inf_{z \in X_{L}}|x-z|<\frac{\pi}{2},$$
    for any $L$ big enough, then
    we can deduce (\ref{metric}) which is
    equivalent to
    $$\sup_{\omega\in \mathbb{S}^{d}}     \min_{z \in L\exp^{-1}(\rho_{L}\mathcal{Z}(L))}|L \exp^{-1}(\omega)-z|<\frac{\pi}{2}.$$
    But this follows from the condition on the mesh norm the property
    $$d(\exp \frac{x}{L},\exp \frac{y}{L})=\frac{|x-y|}{L}+o(L^{-1}),\;\;\;\;x,y\in \R^{d}.$$

    Therefore the condition on $\mathcal{Z}$ implies that it is an $L^\infty$-MZ array.
    In order to deduce the result for $1\le p<\infty$ we define, for $\delta>0,$ the associated
    arrays
    $\mathcal{Z}_{\delta},\mathcal{Z}_{-\delta}$ by
    $\mathcal{Z}_{\delta}(L)=\mathcal{Z}([(1+\delta)L]),\mathcal{Z}_{-\delta}(L)=\mathcal{Z}([(1-\delta)L]).$

    It was proved in \cite[Lemma 4.9.]{Mar07} that if $\mathcal{Z}$ is an $L^\infty$-MZ array then
    $\mathcal{Z}_{\delta}$ is an $L^p$-MZ array for all $1\le p<\infty.$

    Suppose that
    $$\eta=\sup_{L\ge 0} L\rho(\mathcal{Z}(L))<\frac{\pi}{2},$$
    and $\delta>0$ be such that
    $$\eta<\frac{\pi}{2}-\delta\pi.$$

    For $L$ big enough with $L\delta>1$ and any $u\in \mathbb{S}^{d}$ we have
    $$Ld(u,\mathcal{Z}_{L_{1-\delta}})\le \frac{L}{L_{1-\delta}}\eta<\frac{\pi}{2},$$
    therefore $\mathcal{Z}_{-\delta}$ is $L^{\infty}$-MZ and
    $(\mathcal{Z}_{-\delta})_{\delta}=\mathcal{Z}$ is $L^{p}$-MZ.
\qed



{\small
\begin{thebibliography}{S 99}
\bibitem[BCG83]{BCG83} C. Betz, G. A. Cámera, H. Gzyl, {\em Bounds for the first eigenvalue of a spherical cap},
Appl. Math. Optim. 10, no. 3, 193-202, 1983.
\bibitem[Beu89]{Beu89} A. Beurling, {\em The collected works Arne Beurling,} vol. 2, Contemporary Mathematicians,
    Birkhäuser Boston Inc., MA, 1989,
    Complex analysis, Edited by L. Carleson, P. Malliavin, J. Neuberger and J. Wermer.
\bibitem[BC73]{BC73} A. Bonami, J-L. Clerc, {\em Sommes de Cesaro et multiplicateurs
    des developpements en harmoniques spheriques}, Trans. Amer.
    Math. Soc. vol. 183, 223-263, 1973.
\bibitem[DM05]{DM05} S. B. Damelin, V. Maymeskul, {\em On point energies, separation radius and mesh norm for $s$-extremal
    configurations on compact sets in $\R^n,$}
    J. Complexity 21, no. 6, 845-863, 2005.
\bibitem[FM11]{FM11} F. Filbir, H. N. Mhaskar, {\em Marcinkiewicz-Zygmund measures on manifolds}, J. Complexity 27, no. 6, 568-596, 2011. 
\bibitem[FM10]{FM10} F. Filbir, H.N. Mhaskar, {\em A quadrature formula for diffusion polynomials corresponding to a generalized 
heat kernel}, J. Fourier Anal. Appl. 16 (5), 629-657, 2010.
\bibitem[FH76]{FH76} S. Friedland, W. K. Hayman, {\em Eigenvalue inequalities for the Dirichlet problem on spheres and the growth of subharmonic functions,}
Comment. Math. Helv. 51, no. 2, 133-161, 1976.
\bibitem[Ing36]{Ing36} A.E. Ingham, {\em Some trigonometrical inequalities with applications in the theory of series,} Math. Z. 41 (1), 367-379, 1936.
\bibitem[KL05]{KL05} V. Komornik, P. Loreti, {\em Fourier series in control theory},
Springer Monographs in Mathematics. Springer-Verlag, New York, 2005. 
\bibitem[KS98]{KS98} A. B. J. Kuijlaars, E. B. Saff, {\em Asymptotics for minimal discrete energy on the sphere}, Trans. Amer. Math. Soc. 350, no. 2, 523-538, 1998.
\bibitem[MZ37]{MZ37} J. Marcinkiewicz, A. Zygmund, Mean values of trigonometrical
    polynomials, Fund. Math., 28, 131-166, 1937.
\bibitem[Mar07]{Mar07} J. Marzo, {\em Marcinkiewicz-Zygmund inequalities and interpolation by spherical harmonics,}
    J. Funct. Anal. 250, no. 2, 559-87, 2007.
\bibitem[MOC10]{MOC10} J. Marzo, J. Ortega-Cerd\`a, {\em Equidistribution of the Fekete Points on the Sphere,} Constr. Approx. 32, no. 3, 513-521, 2010. 
\bibitem[MNW00]{MNW00} H.N. Mhaskar, F. J. Narcowich, J. D. Ward \emph{Spherical Marcinkiewicz-Zygmund Inequalities and
        positive quadrature}, Math. of Comp. vol. 70, 1113-1130,
        2000.
\bibitem[Mha06]{Mha06} H. N. Mhaskar, {\em Weighted quadrature formulas and approximation by zonal function networks on the sphere,}
    J. Complexity, {\bf 22} (3), 348-370, 2006.
\bibitem[MNSW02]{MNSW02} H.N. Mhaskar, F. J. Narcowich, N. Sivakumar, J. D. Ward, {\em Approximation with interpolatory constraints,}
         Proc. Amer. Math. Soc. , 130, no. 5, 1355-1364, 2002.
\bibitem[Nik75]{Nik75} S. M. Nikol'ski\u{\i}, {\em Approximation of Functions of Several Variables
     and Imbedding Theorems}, Springer-Verlag, Die Grundlehren der mathematische
     Wis\-sen\-schaf\-ten in Einzeldarstellungen, 205, 1975.
\bibitem[OU10]{OU10} A. Olevskii, A. Ulanovskii, {On Ingham-type interpolation in $\mathbb{R}^n$}, C. R. Math. Acad. Sci. Paris 348, no. 13-14, 807-810, 2010. 
\bibitem[OCP11]{OCP11} J. Ortega-Cerd\`a, B. Pridhnani, {\em Carleson measures and Logvinenko-Sereda sets on compact manifolds}, 
Forum Math. 25, Issue 1, 151-172, 2013.
\bibitem[OCP12]{OCP12} J. Ortega-Cerd\`a, B. Pridhnani, {\em  Beurling-Landau's density on compact manifolds}, J. Funct. Anal. 263, 2102-2140, 2012.
\bibitem[Pin81]{Pin81} M. A. Pinsky, {\em The first eigenvalue of a spherical cap,}
Appl. Math. Optim. 7 , no. 2, 137-139, 1981.
\bibitem[Rei90]{Rei90} M. Reimer, {\em Constructive Theory of Multivariate Functions,} BI Wissenschaftsverlag, Mannheim,
    1990.
\bibitem[Ron74]{Ron74} L. I. Ronkin, {\em Introduction to the theory of entire functions of several variables},
    Trans. Math. Monographs, AMS, Providence, RI, 1974.
\bibitem[Sei95]{Sei95} K. Seip, {\em On the connection between exponential bases and certain related sequences in $L^2(-\pi,\pi)$}, 
J. Funct. Anal. 130, 131-160, 1995.
\bibitem[SW04]{SW04} I. H. Sloan, R. S. Womersley, {\em Extremal systems of points and numerical integration on the sphere}, Adv. Comput. Math. 21, no. 1-2, 107-125, 2004.
\bibitem[Sze39]{Sze39} G. Szeg\"o, {\em Orthogonal polynomials}, American Mathematical Society, Colloquium Publications, vol. 23,
1939.
\end{thebibliography}
}
\end{document}